\documentclass{amsart}
\usepackage{amssymb}

\newtheorem{thm}{Theorem}[section]
\newtheorem{cor}[thm]{Corollary}
\newtheorem{lem}[thm]{Lemma}

\theoremstyle{definition}

\theoremstyle{remark}

\numberwithin{equation}{section}

\begin{document}
\title[Absolutely Flat Idempotents]{Absolutely Flat Idempotents}

\author[J.M. Groves]{Jonathan M. Groves}
\address{Department of Mathematics, Austin Peay State University, Clarksville, TN 37044}
\email{jonny77889@yahoo.com}
\author[Y. Harel]{Yonatan Harel}
\address{Institute of Mathematics, The Hebrew University, 91904
Jerusalem, Israel}
\email{zapkif@hotmail.com}
\author[C.J. Hillar]{Christopher J. Hillar}
\address{Department of Mathematics, University of California, Berkeley, CA 94720}
\email{chillar@math.berkeley.edu}
\author[C.R. Johnson]{Charles R. Johnson}
\address{Department of Mathematics, College of William and Mary, Williamsburg, VA
23187-8795} \email{crjohnso@math.wm.edu}
\author[P.X. Rault]{Patrick X. Rault}
\address{Department of Mathematics, College of William and Mary, Williamsburg, VA
23187-8795}
\email{pxraul@wm.edu}

\thanks{This research was conducted during the summer of 2002 at the
College of William and Mary's Research Experiences for
Undergraduates program and was supported by NSF REU grant
DMS-96-19577. The contribution of the second author is part of his
``Amirim'' research project in mathematics, prepared at the Hebrew
University of Jerusalem.  The work of the third author is
supported under a National Science Foundation Graduate Research
Fellowship.}


\begin{abstract}
A real $n$-by-$n$ idempotent matrix $A$ with all entries having
the same absolute value is called {\it absolutely flat}. We
consider the possible ranks of such matrices and herein
characterize the triples: size, constant, and rank for which such
a matrix exists. Possible inequivalent examples of such matrices
are also discussed.
\end{abstract}

\maketitle


\section{Introduction}

We consider a problem, suggested, in part by \cite{Zip1} and
specifically mentioned to the other authors by Harel. The problem
considered in \cite{Zip1} (see also \cite{Zip2,Zip3,Zip4}) is
about the isomorphic classification of the ranges of nicely
bounded projections in some classical Banach spaces. It has been
solved in \cite{Zip1} in the special case of projections of small
norms and another special case is that of absolutely flat
idempotents. We also found this question of independent interest.

For which positive integers $n$, does there exist an $n$-by-$n$
real, idempotent matrix $A$ of rank $r$, all of whose entries are
a positive constant $c$ in absolute value ({\it absolutely flat})?
From the equation $A^{2}=A$, it readily follows that $c$ must be
$1/k$ for some positive integer $k\leq n$. Thus, the key
parameters of our problem are $n,k,r$: for which triples of
positive integers is there a matrix $A$ of desired type? Since
$kA=B$ is a $\pm 1$ matrix, an equivalent formulation concerns the
existence of a $\pm 1$ matrix $B$ such that $B^{2}=kB$, and we
reserve the letter $B$ for such a $\pm 1$ matrix that comes from a
given $A$ in sections 2 and 3 below.

Since the minimal polynomial of $A$ must divide (and, in fact,
equal in our case) $x^{2}-x$, $A$ is diagonalizable (\cite[p.
145]{HJ1}) and all of its eigenvalues must be 0 or 1. Importantly,
Tr[$A$] = rank $A$, as each is simply the count of the number of
eigenvalues equal to 1. We first derive two number theoretic
necessary conditions that constrain feasible triples $n,k,r$.
Then, we show that for odd $n$, only $r=1$ is possible and that
all triples $n,k,1$ meeting the necessary conditions do occur.
Finally, for even $n$, all triples meeting the necessary
conditions occur, completing a characterization of feasible
triples. We also discuss the existence of multiple matrices,
distinct modulo obvious symmetries of the problem, which are
absolutely flat idempotents for the same parameters $n,k,r$.

\section{The Elementary Necessary Conditions}

A signature matrix is a diagonal matrix $S$ with diagonal entries $\pm 1$.
It is clear that similarity does not change the property of idempotence.
Further, permutation and signature similarity do not change the set of
absolute values of the entries of a matrix. Thus, permutation and signature
similarity do not change whether $A$ is an absolutely flat idempotent, nor
do they change the parameters $n,k,r$ if $A$ is.

From the equality of rank and trace, an absolutely flat idempotent
must have at least one positive diagonal entry, and, therefore,
after a permutation similarity, we may assume a positive number as
the (1,1) entry of $A$. Then, any absolutely flat idempotent may
be normalized, by signature similarity, so that all the entries in
its first column are positive. We generally assume this
normalization. From $A^{2}=A$, it then follows that the number of
negative entries in each row is constant. Call this number $u\geq
0$. Similarly, let $m\geq 0$ denote the number of negative entries
on the main diagonal of $A$ (an absolutely flat idempotent of
parameters $n,k,r$). The trace of $A$ is $[(n-m)/k]-[m/k]$ =
$(n-2m)/k$, but since rank equals trace, we have $n-2m=kr$ or
\begin{equation}
n=rk+2m,  \label{elem1}
\end{equation}
the first of our necessary conditions. The second follows from $B^{2}=kB$,
with $B$ in normalized form. The inner product of the first (any) row of $B$
with the first (normalized) column has $n-u$ positive summands and $u$
negative summands, with the net sum being $k$. Thus, $n-2u=k$ or
\begin{equation}
n-k=2u.  \label{elem2}
\end{equation}
Since $m\geq 0$, it follows from (\ref{elem1}) that
\begin{equation}
rk\leq n,  \label{elem3}
\end{equation}
and it follows from (\ref{elem2}) that
\begin{equation}
n\text{ and }k\text{ have the same parity}.  \label{elem4}
\end{equation}
It also follows from (\ref{elem1}) that $n$ is odd if and only if $r$ and $k$
are odd.

\section{The Odd Case}

Many triples with $n$ odd and $r > 1$ (and necessarily odd) satisfy the
requirements (\ref{elem1}) and (\ref{elem2}). However, interestingly,
absolutely flat idempotents never exist in such cases.

\begin{thm}
\label{thm1} For an odd integer $n$, there is an absolutely flat
idempotent $A$ with parameters $n,k,r$ if and only if $r=1$ and
$k\leq n$ is odd. In this event, the matrix $A$ is unique up to
signature/permutation similarity;
\begin{equation*}
A=\frac{1}{k}\left[ {\begin{array}{*{20}c} 1 & \cdots & 1 & { - 1}
& \cdots & { - 1} \\ \vdots & \ddots & \vdots & \vdots & \ddots &
\vdots \\ 1 & \cdots & 1 & { - 1} & \cdots & { - 1} \\
\end{array}}\right] ,
\end{equation*}
in which there are $m=u=(n-rk)/2$ columns of -1's.
\end{thm}

\begin{proof} If $r=1,$ and $k\leq n$ is odd, it is easily
checked that the displayed matrix $A$ shows existence.
Furthermore, in this event, any absolutely flat idempotent that is
normalized via (permutation and) signature similarity to have
positive first column and then by permutation similarity to have
all positive entries in the first $n-u$ columns, will have all
rows equal and appear as the displayed $A$. It follows that
$m=u=(n-rk)/2$ and that this is the number of negative columns.

If $n$ is odd, we already know that $k\leq n$ is odd and $r$ is
odd. We show that $r=1$ in two cases: $k=1;k>1$. Let
$(n,k,r)=(2l+1,2t+1,2s+1).$ Consider $B$ in normalized form, so
that $B^{2}=kB$, with $B$ a $\pm 1$ matrix, and partition $B$ as
\begin{equation}
B=\left[
\begin{array}{cc}
1 & f^{T} \\
e & C
\end{array}
\right]
\end{equation}
in which $e$ is the $2l$-by-$1$ vector of $1$'s and $f$ is a
$2l$-by-$1$ vector consisting of $(l+t)$ $1$'s followed by $(l-t)$
-1's. From $B^{2}=kB$, it follows that
\begin{equation}\label{4system}
\begin{split}
&f^T e = 2t
\\
&f^T C = 2tf^T
\\
&Ce = 2te
\\
&ef^T + C^2 = (2t+1)C.
\end{split}
\end{equation}
Multiplication of both sides of the last equation on the left by $C$ and use
of $Ce=2te$ yields
\begin{equation}\label{minpolyeq}
\begin{split}
0 &= C^3 - (2t+1)C^2 + 2tef^T
\\
&= C^3 - (4t+1)C^2 + 2t(2t+1)C.
\end{split}
\end{equation}
Since $C$ is a $\pm 1$ matrix and is of even dimension, it follows
that $2\mid C^{2}$ (entry-wise) and then, by a simple induction,
that $2^{q} \mid C^{2^{q}}$ for each positive integer $q$. Thus,
$2^{q}\mid$ Tr$\left[C^{2^{q}}\right]$ for all positive integers
$q$.

Now distinguish two possibilities: $k=1$ $(t=0)$; and $k>1$
$(t>0)$. In the former case, (\ref{minpolyeq}) gives $C^{3}=C^{2}$
and, thus, by induction, $C^{2^{q}}=C^{2}.$ Therefore, $2^{q}\mid
C^{2}$ for all positive integers $q$, which gives $C^{2}=0$. But
then $C$ is nilpotent; Tr$[C]=0$, and Tr[$B$] = Tr[$A$] = rank $A$
= $r = 1$, as was to be shown.

Now, suppose $t>0$ $(k>1)$. First, rank $C$ = rank $B$, as the
first column of $B$ is $1/2t$ times the sum of the last $2l$
columns of $B$ (by the first and third equations of
(\ref{4system})), and the equation, $f^TC = 2tf^T$, implies that
$f^T$ can be written as a linear combination of rows of $C$. From
(\ref{minpolyeq}), $C$ is diagonalizable with distinct eigenvalues
from $\left\{0,2t,2t+1\right\} $. \ Let $a\geq 0$ be the number of
eigenvalues of $C$ equal to $2t$, $b\geq 0$ be the number equal to
($2t+1);$ then there are $2l-a-b$ of them equal to 0. Since rank
$C$ = rank $B$, we have
\begin{equation*}
a+b=2s+1.
\end{equation*}
Also, $r=$ Tr[$B]/(2t+1)$, so that $2ta+(2t+1)b =$ Tr[$C$] =
Tr[$B]-1=(2s+1)(2t+1)-1$, or $2ta+(2t+1)b=(2s+1)2t+2s$. These two
equations have the unique solution $(a,b)=(1,2s)$.  We may now
calculate Tr$\left[C^{2^{q}}\right]$ as
\begin{equation*}
(2t)^{2^{q}}+2s(2t+1)^{2^q}.
\end{equation*}
Since $2^{q}\mid C^{2^{q}}$, still, and thus $2^{q}\mid$
 Tr$\left[C^{2^{q}}\right]$ for all positive integers $q$, we have $s=0$,
or $r=1,$ as was to be shown. This concludes the proof.
\end{proof}

\section{The Even Case}

When $n$ is even, conditions (\ref{elem1}) and (\ref{elem2}) still govern
existence, but the overall situation is remarkably different from the odd
case. Now, there is existence whenever the conditions are met. Here we
exhibit an absolutely flat idempotent for each triple $n,k,r$ meeting the
conditions (2.1) and (2.2).

For a given positive integer $k$, define
\begin{equation}  \label{defPM}
P = \frac{1}{{k}}\left[ {\begin{array}{*{20}c} 1 & 1 \\ 1 & 1 \\
\end{array}}
\right], \ \ \ M = \frac{1}{{k}}\left[ {\begin{array}{*{20}c} 1 & { - 1} \\
1 & { - 1} \\ \end{array}} \right].
\end{equation}
Because of (\ref{elem2}), $k = 2t$ must be even, and we have
\begin{equation}  \label{PMfact}
P^2 = \frac{1}{t}P, \ PM = \frac{1}{t}M, \ M^2 = 0, \ MP = 0.
\end{equation}
Solutions may now be constructed using the $P$'s and $M$'s as blocks. For
example, a solution for $(n,k,r)$ = (8,2,3) is
\begin{equation*}
\left[ {\begin{array}{*{20}c} P & M & M & M \\ M & P & M & M \\ M
& M & P & M \\ P & M & M & M \\ \end{array}} \right].
\end{equation*}
Since $n$ is even and, therefore, $k$ is even, we assume our
parameters are of the form $(n,k,r)$ = $(2l,2t,r)$; $r$ need not
be even. From (\ref{elem3}), it follows that $tr \leq l$. As proof
of the following theorem, we give a general strategy for
constructing absolutely flat idempotents with parameters $2l,2t$
and $r$, $tr \leq l$.

\begin{thm}
\label{thm2} Let $n = 2l$, $k = 2t$ and $r$ be positive integers
such that $tr \leq l$. Then, there is an absolutely flat
idempotent with parameters $n,k,r$. In particular, whenever $n$ is
even, there is an absolutely flat idempotent whenever conditions
(\ref{elem1}) and (\ref{elem2}) are met.
\end{thm}

\begin{proof}
Let $(n,k,r) = (2l,2t,r)$.  By the elementary necessary
conditions, express $2l$ as $2tr + 2m$ for some $m \in \mathbb N.$
Let $P$ and $M$ be the matrices as in (\ref{defPM}).  Examine now
the block matrix,
\begin{equation}
A = \left[ {\begin{array}{*{20}c}
   P &  \cdots  & P & {} & {} & {} & {} & {} & {} & {}  \\
    \vdots  &  \ddots  &  \vdots  & {} & {} & {} & {} & {} & {} & {}  \\
   P &  \cdots  & P & {} & {} & {} & {} & {} & {} & {}  \\
   {} & {} & {} &  \ddots  & {} & {} & {} & {} & {} & {}  \\
   {} & {} & {} & {} & P &  \cdots  & P & {} & {} & {}  \\
   {} & {} & {} & {} &  \vdots  &  \ddots  &  \vdots  & {} & {} & {}  \\
   {} & {} & {} & {} & P &  \cdots  & P & {} & {} & {}  \\
   P &  \cdots  & P & {} & {} & {} & {} & M &  \cdots  & M  \\
    \vdots  &  \ddots  &  \vdots  & {} & {} & {} & {} &  \vdots  &  \ddots
&  \vdots   \\
   P &  \cdots  & P & {} & {} & {} & {} & M &  \cdots  & M  \\
\end{array}} \right].
\end{equation}

The matrix $A$ consists of $r$ $t$-by-$t$ blocks of $P$'s along
the main diagonal and an $m$-by-$t$ block of $P$'s in the lower
left-hand corner.  All other blocks in $A$ are $M$'s.  It is then
an elementary exercise in block matrix multiplication (using
(\ref{PMfact})) that $A^2 = A$. As the trace of $A$ is $r$, it
follows that the rank of $A$ is $r$. This completes the proof.
\end{proof}

Theorems \ref{thm1} and \ref{thm2} provide a complete characterization of
the triples $n,k,r$ for which absolutely flat idempotents exist. We note
that any positive integer $k$ ($n$) may occur, but for odd $k$ ($n$), only
rank 1 matrices exist. On the other hand, any rank may occur. In either
case, $n$ need be sufficiently large.

\section{Multiple Solutions}

By appealing to the Jordan canonical form, any two $n$-by-$n$ idempotents of
the same rank are similar. However, for our problem, restriction to
permutation and signature similarity is more natural; of course, permutation
and signature similarities send one solution for $n,k,r$ to another for the
same $n,k,r$. Although it has not been important for our earlier results,
transposition is another natural operation sending one solution to another.
It is natural to ask how many solutions, distinct up to permutation,
signature similarity, and transposition can occur. When $r=1$, it is easily
worked out that there is only one (when there is one). The form mentioned in
Theorem \ref{thm1} is canonical (even when $n$ is even).

However, already for the parameters (8,2,2), there can be distinct
solutions. For example,
\begin{equation*}
A_{1}=\frac{1}{2}\left[ {\begin{array}{*{20}c} 1 & 1 & 1 & 1 & 1 &
{ - 1} & {-1} & { - 1} \\ 1 & 1 & {-1} & { - 1} & {-1} & 1 & 1 & 1
\\ 1 & 1 & 1 & 1 & 1 &
{ - 1} & {-1} & { - 1} \\ 1 & 1 & 1 & 1 & 1 & { - 1} & {-1} & { -
1} \\ 1 & 1 & 1 & 1 & 1 & { - 1} & {-1} & { - 1} \\ 1 & 1 & 1 & 1
& 1 & { - 1} & {-1} & { - 1} \\
1 & 1 & 1 & 1 & 1 & { - 1} & {-1} & { - 1}
\\ 1 & 1 & {-1} & { - 1} & {-1} & 1 & 1 & 1 \\ \end{array}}\right]
\end{equation*}
and
\begin{equation*}
A_{2}=\frac{1}{2}\left[ {\begin{array}{*{20}c} 1 & 1 & 1 & {-1} &
1 & 1 & {-1} & { - 1} \\ 1 & 1 & 1 & {-1} & {-1} & { - 1} & 1 & 1
\\ 1 & 1 & 1 & {-1} & {-1} & { - 1} & 1 & 1 \\
1 & 1 & 1 & {-1} & {-1} & { - 1} & 1 & 1
\\ 1 & 1 & 1 & {-1} &
1 & 1 & {-1} & { - 1} \\ 1 & 1 & 1 & {-1} & {-1} & { - 1} & 1 & 1
\\ 1 & 1 & 1 & {-1} & {-1} & { - 1} & 1 & 1 \\
1 & 1 & 1 & {-1} & {-1} & { - 1} & 1 & 1 \\ \end{array}}\right]
\end{equation*}
are both absolutely flat (8,2,2) idempotents. To see that $A_{1}$
is not permutation/signature similar to $A_{2}$ (or its
transpose), we mention an idea that we used to discover some of
the construction herein, but was not needed in the proofs thus
far. We say that two rows (columns) of an $n$-by-$n $ $\pm 1/k$
matrix are {\it of the same type} if they are either identical or
negatives of each other. It is an easy exercise that the number of
distinct row types (number of distinct column types) is unchanged
by either signature similarity or permutation similarity.

Additionally, we define the {\it row} ({\it column}) {\it
multiplicity} of an absolutely flat idempotent matrix, $A$, to be
the multiset consisting of the number of rows (columns) for each
row (column) type. It is again an easy exercise that
permutation/signature similarity does not change the row (column)
multiplicity of an absolutely flat idempotent matrix.  In the
matrix $A_{1}$ the row (column) multiplicity is $\{6,2\}$
($\{6,2\}$), while in $A_{2}$ the row (column) multiplicity is
$\{6,2\}$ ($\{4,4\}$). Thus, $A_{1}$ cannot be transformed to
$A_{2}$ by any combination of permutation/signature similarities
and/or transposition (though they are similar).

Of course, the number of row types in a matrix is at least the
rank. We note that the construction technique of Theorem
\ref{thm2} always produces a solution with the same number of row
types as rank. The (8,2,3) example, $A_3$, below demonstrates that
larger numbers of row types are possible. However, it may be shown
that for rank 2 absolutely flat idempotents, only 2 row and column
types are possible.

\begin{equation*}
A_{3}=\frac{1}{2}\left[ {\begin{array}{*{20}c} 1 & 1 & 1 & 1 & 1 &
{ - 1} & {-1} & { - 1} \\ 1 & 1 & 1 & 1 & 1 & { - 1} & {-1} & { -
1} \\ 1 & 1 & 1 & 1 & 1 & { - 1} & {-1} & { - 1} \\ 1 & 1 & 1 & 1
& 1 & { - 1} & {-1} & { - 1}
\\ 1 & 1 & {-1} & {-1} & 1 & 1 & 1 & { - 1} \\ 1 & 1 & {-1} & 1 & 1 & 1 &
{-1} & { - 1} \\ 1 & 1 & 1 & { - 1} & 1 & { - 1} & 1 & { - 1} \\ 1
& 1 & 1 & 1 & 1 & { - 1} & {-1} & { - 1} \\ \end{array}}\right]
\end{equation*}

\begin{lem}
\label{rank2rows} A rank 2 absolutely flat matrix has precisely 2 row types
and 2 column types.
\end{lem}

\begin{proof}
We prove the result for row types as the case of columns is
similar. Let $A$ be a rank 2 absolutely flat matrix. Performing
permutation and signature similarity we may assume that the first
column of $A$ consists only of positive entries, as this doesn't
change the number of row types.  Since $A$ has rank 2, there are
at least 2 distinct row types.  Let $x$ and $y$ be the two rows
corresponding to these row types, and let $w$ be an arbitrary
other row in $A$.  Then,
\[w = ax + by\] for some $a,b \in \mathbb Q$. Clearly, we must
have $a+b = 1$ because the initial entries of $w$, $x$, and $y$
are all the same. Since $x$ and $y$ are different rows, it follows
from the absolutely flat property that $a-b = 1$ or $a-b = -1$. In
the first case, we have $a = 1$ and $b=0$, and in the second, it
follows that $a = 0$ and $b=1$. This completes the proof.
\end{proof}

We now consider the problem of counting all different rank 2
absolutely flat idempotent matrices. As we are interested in
distinct solutions up to permutation and signature similarity, we
first put our matrix in a normalized form. Let $A$ be a rank 2
absolutely flat idempotent matrix with parameters $(n,k,2)$. As
before, we can perform a permutation and signature similarity to
make the first column of $A$ positive. Let $a$ ($b$) be the number
of all positive (negative) columns of $A$. Through another
permutation similarity, we may assume that the first $a$ columns
of $A$ are positive and that the next $b$ columns of $A$ are
negative. From Lemma \ref {rank2rows}, the remaining $n-a-b$
columns of $A$ are of one type. Let $v$ be one of these columns
(necessarily containing both a positive and a negative entry) and
let $c$ be the number of them in $A$. Notice that the other $d =
n-a-b-c$ columns must be $-v$. Since $c + b = u$ and $d+ b = u$
(the number of negative entries in each row must be $u$), it
follows that $d = c$. These normalizations partition our matrix as
\begin{equation}  \label{normalform}
F = \left[ {\begin{array}{*{20}c} {P_{a,a} } & {M_{a,b} } &
{W_{a,c} } & { -
W_{a,c} } \\ {P_{b,a} } & {M_{b,b} } & {X_{b,c} } & { - X_{b,c} } \\
{P_{c,a} } & {M_{c,b} } & {Y_{c,c} } & { - Y_{c,c} } \\ {P_{c,a} }
& {M_{c,b} } & {Z_{c,c} } & { - Z_{c,c} } \\ \end{array}} \right].
\end{equation}

Here, the $P_{i,j}$ are positive matrices of sizes $i$-by-$j$;
$M_{i,j}$ are negative matrices of sizes $i$-by-$j$; and
$W_{i,j}$, $X_{i,j}$, $Y_{i,j}$, and $Z_{i,j}$ are matrices of
sizes $i$-by-$j$ with exactly 1 row type. Through further
permutation, it is clear that the columns of $W_{a,c}$, $X_{b,c}$,
$Y_{c,c}$, and $Z_{c,c}$ can be assumed to begin with all positive
entries and end with all negative ones:
\begin{equation}  \label{signform}
\left[ {\begin{array}{*{20}c} + & \cdots & + \\ \vdots & \ddots & \vdots \\
+ & \cdots & + \\ - & \cdots & - \\ \vdots & \ddots & \vdots \\ -
& \cdots & - \\ \end{array}} \right].
\end{equation}
Let $a_p$ ($a_m$) be the number of positive (negative) rows in
$W_{a,c}$; $b_p$ ($b_m$) be the number of positive (negative) rows
in $X_{b,c}$; $c_{1p}$ ($c_{1m}$) be the number of positive
(negative) rows in $Y_{c,c}$; and $c_{2p}$ ($c_{2m}$) be the
number of positive (negative) rows in $Z_{c,c}$. The final matrix
produced after this sequence of operations is called the {\it
standard} form of $A$.

We now derive necessary conditions on the parameters defined above
for the matrix as in (\ref{normalform}) to be idempotent. Clearly,
we must have $a_p+a_m = a$, $b_p+b_m = b$, $c_{1p}+c_{1m} = c$,
$c_{2p}+c_{2m} = c$, $a+b+2c = n$, and $b+c = u$. Examining the
inner product of the first row and the first column, we see that
$a-b = k$, and looking at the inner products of each row type with
the second column type produces the equations,
\begin{equation}  \label{innerprod1}
a_p -a_m-b_p+b_m+c_{1p}-c_{1m}-c_{2p}+c_{2m} = k
\end{equation}
and
\begin{equation}  \label{innerprod2}
a_p -a_m-b_p+b_m-c_{1p}+c_{1m}+c_{2p}-c_{2m} = -k.
\end{equation}
Adding equations (\ref{innerprod1}) and (\ref{innerprod2}) gives us that
\begin{equation*}
\begin{split}
0 &= a_p/2-a_m/2+b_m/2 -b_p/2 \\
&= a_p - b_p + (b-a)/2 \\
&= a_p - b_p - k/2,
\end{split}
\end{equation*}
and a similar computation with the subtraction of
(\ref{innerprod1}) and (\ref{innerprod2}) produces the equation,
$c_{1p} = c_{2p} + k/2$. Many of these necessary conditions are
actually redundant, and so we will only consider the system,
\begin{equation}  \label{necsystem}
\begin{split}
&a - b = k \\
&b+c = u \\
&a_p = b_p + k/2 \\
&c_{1p} = c_{2p} + k/2.
\end{split}
\end{equation}
In fact, we have the following

\begin{thm}
\label{sufthm} A matrix in standard form in which
$a,b,c,a_p,b_p,c_{1p},c_{2p}$ are all nonnegative and satisfy
(\ref{necsystem}) is an (n,k,2) idempotent.
\end{thm}

\begin{proof}
Assume that $A$ is in standard form with
$a,b,c,a_p,b_p,c_{1p},c_{2p} \geq 0$ and (\ref{necsystem})
satisfied.  To prove idempotence, we need to check three inner
products.  The inner product of the first row type and first
column type is just $a-b = k$, and the inner product of the first
row type and the second column type is
\[a_p - (a-a_p) -b_p +(b-b_p) +(c_{2p}+k/2) -(c-c_{2p}-k/2) -
c_{2p} + (c-c_{2p})\]
\begin{equation*}
\begin{split}
&= 2a_p-2b_p-a+b+k
\\
&= k
\end{split}
\end{equation*}
as desired.  A similar computation involving the second row type
and the second column type gives us
\[a_p -(a-a_p) - b_p +(b-b_p) - (c_{2p}+k/2) + (c-c_{2p}-k/2) + c_{2p}
-(c-c_{2p})\]
\begin{equation*}
\begin{split}
&= 2a_p - 2b_p + b -a - k
\\
&= -k.
\end{split}
\end{equation*}
Finally, adding the equations $2b+2c = 2u$ and $a-b = k$ gives us
that $a+b+2c = n$, completing the proof.
\end{proof}

In what follows, the multiplicities of an absolutely flat
idempotent matrix will be important.  Let $x_A =
a_p+b_p+c_{1p}+c_{2p} = 2b_p+2c_{2p} + k$ and set $y_A = a+b$.
Then, the row and column multiplicities of $A$ in standard form
are $\{x_A,n-x_A\}$ and $\{y_A,n-y_A\}$, respectively. The
following lemma is a natural consequence of the symmetries of the
problem.

\begin{lem}
\label{swaplemma} Let $A$ be a rank 2 absolutely flat idempotent
in standard form with row and column multiplicities of
$\{x_A,n-x_A\}$ and $\{y_A,n-y_A\} $ as above. Then, $A$ is
permutation/signature equivalent to a matrix $B$ in standard form
with $x_B = n-x_A$ and $y_B = y_A$. Similarly, $A $ is
permutation/signature equivalent to a matrix $B$ in standard form
with $y_B = n-y_A$ and $x_B = x_A$.
\end{lem}

\begin{proof}
Let $A$ be as in (\ref{normalform}).  After permuting the last
$2c$ columns and the corresponding last $2c$ rows, $A$ becomes
\[\left[ {\begin{array}{*{20}c}
   {P_{a,a} } & {M_{a,b} } & {-W_{a,c} } & { W_{a,c} }  \\
   {P_{b,a} } & {M_{b,b} } & {-X_{b,c} } & { X_{b,c} }  \\
   {P_{c,a} } & {M_{c,b} } & {-Z_{c,c} } & { Z_{c,c} }  \\
   {P_{c,a} } & {M_{c,b} } & {-Y_{c,c} } & { Y_{c,c} }  \\
   \end{array}} \right].\]  Through further permutation,
the columns of $-W_{a,c}$, $-X_{b,c}$, $-Z_{c,c}$, and $-Y_{c,c}$
can be made to look like those in (\ref{signform}).  Now, this
final matrix, $B$, is in normal form with $x_B = n-x_A$, and $y_B
= y_A$ as desired.

As for the second statement in the lemma, first perform a
signature similarity on $A$ that makes each column of, \[ \left[
{\begin{array}{*{20}c}
   { W_{a,c} }  \\
   { X_{b,c} }  \\
   { Y_{c,c} }  \\
   { Z_{c,c} }  \\
   \end{array}} \right], \] either all positive or all negative,
and then perform a permutation similarity to bring our matrix back
into standard form.  It is clear that this new matrix, $B$, has
$y_B = n-y_A$.  If $x_B = x_A$, then we are done.  Otherwise, $x_B
= n-x_A$, and we can proceed as above to form an equivalent
matrix, $B'$, with $x_{B'} = n-x_{B} = x_A$ and $y_{B'} = y_B =
n-y_A$. This completes the proof of the lemma.
\end{proof}

We are now in a position to give bounds for the number of rank 2 absolutely
flat idempotent matrices up to permutation/signature similarity and
transposition. A straightforward verification (using (\ref{elem1}) and (\ref
{elem2})) shows that
\begin{gather*}
\begin{split}
&a_p = k/2, \ a_m = \left\lceil {n/4} \right\rceil \\
&b_p = 0, \ b_m = \left\lceil {m/2} \right\rceil \\
&c_{1p} = k/2, \ c_{1m} = \left\lfloor {m/2} \right\rfloor \\
&c_{2p} = 0, \ c_{2m} = \left\lfloor {n/4} \right\rfloor
\end{split}
\end{gather*}
satisfy (\ref{necsystem}) and, therefore, produce an $(n,k,2)$ absolutely
flat idempotent matrix by Theorem \ref{sufthm} (this is, in fact, the
solution found in Theorem \ref{thm2}). In this case, the row multiplicity is
$\{k,n-k\}$ and the column multiplicity is $\{2\left\lceil {n/4}
\right\rceil,2\left\lfloor {n/4} \right\rfloor\}$.

Let $a,b,c,a_p,b_p,c_{1p},c_{2p}$ be an arbitrary solution to
(\ref {necsystem}). Set $t = a_p-k/2 = b_p$, $q = c_{1p} -k/2 =
c_{2p}$, and let $l = \left\lceil {n/4} \right\rceil - a_m$. Since
$a_p+a_m-b_p-b_m = k$, it follows that $l = \left\lceil {m/2}
\right\rceil - b_m$. If we set $p = \left\lfloor {m/2}
\right\rfloor - c_{1m}$ and $y = \left\lfloor {n/4} \right\rfloor
- c_{2m}$, then from $c_{1p} + c_{1m} = c_{2p} + c_{2m}$ we must
have $p=y$. Finally, the equation $b+c = u$ implies that $p =
q+t-l$. It is easily seen that these conditions are also
sufficient, and so we have the following.

\begin{thm}
\label{allsoltns} All solutions to (\ref{necsystem}) in nonnegative integers
are given by
\begin{gather*}
\begin{split}
&a_p = k/2 + t, \ a_m = \left\lceil {n/4} \right\rceil - l \\
&b_p = t, \ b_m = \left\lceil {m/2} \right\rceil-l \\
&c_{1p} = k/2+ q, \ c_{1m} = \left\lfloor {m/2} \right\rfloor-q-t+l \\
&c_{2p} = q, \ c_{2m} = \left\lfloor {n/4} \right\rfloor-q-t+l
\end{split}
\end{gather*}
in which $t,q \in \mathbb N$, $l \in \mathbb Z$, and
\begin{equation*}
q+t-\left\lfloor {m/2}\right\rfloor \ \leq \ l \ \leq \ \left\lceil {m/2}
\right\rceil.
\end{equation*}
\end{thm}

In particular, when $m = 0$, we must have $q = t = l = 0$, giving us the
immediate

\begin{cor}
\label{unqsol1} Up to permutation/signature similarity and transposition,
there is only one rank 2 absolutely flat idempotent matrix with $n = 2k$.
\end{cor}

With a careful consideration of Theorem \ref{allsoltns}, we can
produce bounds for the number of inequivalent $(n,k,2)$ absolutely
flat idempotents. Notice that for the parameterized solutions in
Theorem \ref {allsoltns}, we have $x_A = 2t+2q+k$ and $y_A =
k+2\left\lceil {m/2} \right\rceil +2t-2l$.  In particular, the
conditions in Theorem \ref{allsoltns} imply that $x_A = k+2i$ and
$y_A = k+2j$ for some $i,j \in \{0,\ldots,m\}$.

In fact, the converse is true.  Namely, let $i,j \in
\{0,\ldots,m\}$; then, we claim that (\ref{necsystem}) has a
solution, $A$, in which $x_A = k+2i$ and $y_A = k+2j$.  To see
this, fix $i \in \{0,\ldots,m\}$, and let $l$ and $t$ be such that
$0 \leq t \leq i$ and $i - \left\lfloor {m/2}\right\rfloor \leq l
\leq \left\lceil {m/2} \right\rceil$.  Next, set $q = i - t$.
Then, $q,t,l$ gives rise to a solution of (\ref{necsystem}) by
Theorem \ref{allsoltns}, and we have $x_A = k+2i$. Moreover, it is
clear that the value of $t-l$ may be taken to be any number from
$\{-\left\lceil {m/2} \right\rceil,\ldots,\left\lfloor {m/2}
\right\rfloor\}$. This proves the claim.

Since we are looking for inequivalent solutions, we will only
consider (by Lemma \ref{swaplemma}) $i,j \in
\{0,\ldots,\left\lfloor {m/2}\right\rfloor\}$. As transposition
(which switches the row and column multiplicities) could make two
solutions permutation/signature equivalent, it follows from the
discussion above that we have at least
\begin{equation}  \label{lowerbnd}
\sum\nolimits_{j = 0}^{\left\lfloor {m/2} \right\rfloor }
{\sum\nolimits_{i = j}^{\left\lfloor {m/2} \right\rfloor } {{1} }
} = {\binom{\left\lfloor {m/2}\right\rfloor + 2 }{2}}
\end{equation}
inequivalent $(n,k,2)$ absolutely flat idempotents.

We now discuss bounding the number of solutions from above. Given
$i,j \in \{0,\ldots,\left\lfloor {m/2}\right\rfloor\}$ with $i
\geq j$ (recall that transposition may be used to swap row and
column multiplicities), we will count the number of triples,
($q$,$t$,$l$), that give rise to a rank 2 absolutely flat
idempotent, $A$, with $x_A = k+2i$ and $y_A = k+2j$. From Theorem
\ref {allsoltns}, it follows that $t+q = i$ and $\left\lceil {m/2}
\right\rceil+ t-l = j$ in which $t,q \in \mathbb N$ and $i -
\left\lfloor {m/2}\right\rfloor \leq l \leq \left\lceil {m/2}
\right\rceil$. When $l = \left\lceil {m/2} \right\rceil$, we must
have $t = j$ and $q = i-j$, and when $l = \left\lceil {m/2}
\right\rceil - j$, it follows that $t = 0$ and $q = i$.  It is
easy to see, therefore, that there are $j+1$ solutions to such a
system given $i \geq j$. Hence, the total number of inequivalent
solutions is bounded above by,
\begin{equation*}
\sum\nolimits_{j = 0}^{\left\lfloor {m/2} \right\rfloor } {
\sum\nolimits_{i = j}^{\left\lfloor {m/2} \right\rfloor } {\left(
{j + 1} \right)} } = \frac{{\left\lfloor {m/2} \right\rfloor
\left( {\left\lfloor {m/2} \right\rfloor + 1} \right)\left(
{\left\lfloor {m/2} \right\rfloor + 2} \right)}} {6} +
{\binom{\left\lfloor {m/2}\right\rfloor + 2 }{2}}.
\end{equation*}
Combining this computation with (\ref{lowerbnd}) gives us the
following.

\begin{thm}
Let $N$ be the number of inequivalent $(n,k,2)$ absolutely flat idempotent
matrices. Then,
\begin{equation*}
{\binom{\left\lfloor {m/2}\right\rfloor + 2 }{2}} \leq N \leq
\frac{{ \left\lfloor {m/2} \right\rfloor \left( {\left\lfloor
{m/2} \right\rfloor + 1 } \right)\left( {\left\lfloor {m/2}
\right\rfloor + 2} \right)}} {6} + {\binom{\left\lfloor
{m/2}\right\rfloor + 2 }{2}}.
\end{equation*}
\end{thm}

When $m = 1$, it is clear that $N=1$, and thus we have

\begin{cor}
\label{unqsol2} Up to permutation/signature similarity and transposition,
there is only one rank 2 absolutely flat idempotent matrix with $n = 2k+2$.
\end{cor}


\end{document}